\documentclass[a4paper,12pt]{article}
\usepackage[T2A]{fontenc}
\usepackage[english]{babel}
\usepackage{graphicx}
\usepackage{amsmath}
\usepackage{mathrsfs}
\begin{document}
\begin{center}

%%%%%%%%%%%%%%
%% Øàïêà ñòàòüè
%%%%%%%%%%%%%%

%\begin{flushleft}
%{\small {\it ÓÄÊ 517.954}}
%\end{flushleft}

\vskip 5mm {\Large \bfseries A Priori Estimates for Solutions of
Boundary Value Problems for Fractional-Order Equations}

\vskip 5mm {\bf A. A. Alikhanov}

{Kabardino-Balkarian State University, ul. Chernyshevskogo 175,
 Nalchik,  360004,   Russia}

{Institute for Applied Mathematics Research KBSC RAS, ul. Shortanova
89-A, Nalchik,  360000,   Russia}

{\it e-mail: aaalikhanov@gmail.com}

\end{center}
%%%%%%%%%%%%%%%%%%%%%%%%%%%%%%%%%%%%%%%%%%%%%%%%%%%%%%%%%%%%%%%%%%%%%%%

\sloppy

\begin{abstract}
We consider boundary value problems of the first
and third kind for the diffusionwave equation. By using the method
of energy inequalities, we find a priori estimates for the solutions
of these boundary value problems.
\end{abstract}

Fractional calculus is used for the description of a large class of
physical and chemical processes that occur in media with fractal
geometry as well as in the mathematical modeling of economic and
social-biological phenomena \cite[Chap. 5]{Nakhush}. It was proved
in \cite{Nakhush} that fractional differentiation is a positive
operator; this result permits one to obtain a priori estimates for
solutions of a wide class of boundary value problems for equations
with fractional derivatives. The paper \cite{Samko} deals with the
generalization of the differentiation and integration operations
from integer to fractional, real, and complex orders and with
applications of fractional integration and differentiation to
integral and differential equations and in function theory. In the
paper \cite{Shkhan}, an a priori estimate in terms of a fractional
Riemann-Liouville integral of the solution was obtained for the
solution of the first initial-boundary value problem for the
fractional diffusion equation. The general fractional diffusion
equation $(0<\alpha\leq1)$ with regularized fractional derivative
was considered in  \cite{Kochu}.  A more detailed bibliography on
fractional partial differential equations, including the
diffusion-wave equation, can be found, for example, in \cite{Pskhu}.

In the present paper, we use the method of energy inequalities to
obtain a priori estimates for solutions of boundary value problems
for the diffusion-wave equation with Caputo fractional
derivative\cite{Caputo}.

%%%%%%%%%%%%%%%%%%%%%%%%%%%%%%%%%%%%%%%%%%%%%%%%%%%%%%%%%%%%%%%%%%%%%%%%
%%%%%%%%%%%%%%%%%%%%%%%%%%%%%%%%%%%%%%%%%%%%%%%%%%%%%%%%%%%%%%%%%%%%%%%%
\section{BOUNDARY VALUE PROBLEMS FOR THE FRACTIONAL DIFFUSION EQUATION}
%%%%%%%%%%%%%%%%%%%%%%%%%%%%%%%%%%%%%%%%%%%%%%%%%%%%%%%%%%%%%%%%%%%%%%%%
%%%%%%%%%%%%%%%%%%%%%%%%%%%%%%%%%%%%%%%%%%%%%%%%%%%%%%%%%%%%%%%%%%%%%%%%

{\bf 1.1. First Boundary Value Problem.}  In the rectangle $\bar
Q_T=\{(x,t):0\leq x\leq l, 0\leq t\leq T\}$, consider the first
boundary value problem

\begin{equation}\label{ur1}
\partial_{0t}^\alpha u=\frac{\partial}{\partial
x}\left(k(x,t)\frac{\partial u}{\partial x}\right)-q(x,t)u+f(x,t)
,\quad 0<x<l,\quad 0<t\leq T,
\end{equation}

\begin{equation}
u(0,t)=0,\quad u(l,t)=0,\quad 0\leq t\leq T, \label{ur2}
\end{equation}

\begin{equation}
u(x,0)=u_0(x),\quad 0\leq x\leq l, \label{ur3}
\end{equation}
where $$\partial_{0t}^\alpha
u(x,t)=\frac{1}{\Gamma(1-\alpha)}\int_0^t\frac{u_{\tau}(x,\tau)}{(t-\tau)^{\alpha}}d\tau$$
is the Caputo fractional derivative of order $\alpha$ ,
$0<\alpha<1$.

Throughout the following, we assume that there exists a solution
$u(x,t)\in C^{2,1}(\bar Q_T)$ of problem (\ref{ur1})-(\ref{ur3}),
where $C^{m,n}(\bar Q_T)$ is the class of functions that, together
with their partial derivatives of order $m$ with respect to $x$ and
order $n$ with respect to $t$, are continuous on $\bar Q_T$.

Let us prove the following assertion.

{\bf Lemma 1.} For any function $v(t)$ absolutely continuous on
$[0,T]$, one has the inequality
\begin{equation}
v(t)\partial_{0t}^\alpha v(t)\geq \frac{1}{2}\partial_{0t}^\alpha
v^2(t), \quad 0<\alpha<1.
 \label{ur5}
\end{equation}

{\bf Proof.} Let us rewrite inequality (\ref{ur5})  in the form

$$
v(t)\partial_{0t}^\alpha v(t)-\frac{1}{2}\partial_{0t}^\alpha
v^2(t)=\frac{1}{\Gamma(1-\alpha)}v(t)\int\limits_0^t\frac{v_\tau(\tau)d\tau}{(t-\tau)^\alpha}-
\frac{1}{2\Gamma(1-\alpha)}\int\limits_0^t\frac{2v(\tau)v_\tau(\tau)d\tau}{(t-\tau)^\alpha}=
$$

$$
=\frac{1}{\Gamma(1-\alpha)}\int\limits_0^t\frac{v_\tau(\tau)(v(t)-v(\tau))d\tau}{(t-\tau)^\alpha}=
\frac{1}{\Gamma(1-\alpha)}\int\limits_0^t\frac{v_\tau(\tau)d\tau}{(t-\tau)^\alpha}\int\limits_\tau^tv_\eta(\eta)
d\eta=
$$

$$
=\frac{1}{\Gamma(1-\alpha)}\int\limits_0^tv_\eta(\eta)
d\eta\int\limits_0^\eta\frac{v_\tau(\tau)d\tau}{(t-\tau)^\alpha}\equiv
I\geq 0.
$$

Therefore, to prove the lemma, it suffices to show that the integral
$I$ is nonnegative. The integral $I$ takes nonnegative values, since

$$
I=\frac{1}{\Gamma(1-\alpha)}\int\limits_0^{t}v_\eta(\eta)
d\eta\int\limits_0^\eta\frac{v_\tau(\tau)d\tau}{(t-\tau)^\alpha}=
\frac{1}{\Gamma(1-\alpha)}\int\limits_0^{t}(t-\eta)^\alpha\frac{v_\eta(\eta)
d\eta}{(t-\eta)^\alpha}\int\limits_0^\eta\frac{v_\tau(\tau)d\tau}{(t-\tau)^\alpha}=
$$

$$
=\frac{1}{2\Gamma(1-\alpha)}\int\limits_0^{t}(t-\eta)^\alpha\frac{\partial}{\partial\eta}
\left(\int\limits_0^\eta\frac{v_\tau(\tau)d\tau}{(t-\tau)^\alpha}\right)^2d\eta=
$$

$$
=
\frac{\alpha}{2\Gamma(1-\alpha)}\int\limits_0^{t}(t-\eta)^{\alpha-1}
\left(\int\limits_0^\eta\frac{v_\tau(\tau)d\tau}{(t-\tau)^\alpha}\right)^2d\eta\geq
0.
$$

The proof of the lemma is complete.

We use the following notation:  $$\|u\|_0^2=\int_0^l
u^2(x,t)dx,\quad
D_{0t}^{-\alpha}u(x,t)=\frac{1}{\Gamma(\alpha)}\int_0^t\frac{u(x,\tau)}{(t-\tau)^{1-\alpha}}d\tau$$
is the fractional Riemann-Liouville integral of order $\alpha$.

 {\bf Theorem 1.} If $k(x,t)\in C^{1,0}(\bar Q_T)$, $q(x,t)$,
$f(x,t)\in C(\bar Q_T)$, $k(x,t)\geq c_1>0$ and $q(x,t)\geq 0$
everywhere on $\bar Q_T$, then the solution $u(x,t)$ of problem
(\ref{ur1})-(\ref{ur3}) satisfies the a priori estimate

\begin{equation}
\|u\|_0^2+D_{0t}^{-\alpha}\|u_x\|_0^2\leq M
\left(D_{0t}^{-\alpha}\|f\|_0^2+\|u_0(x)\|_0^2\right).
 \label{ur9.1}
\end{equation}

{\bf Proof.} We multiply Eq. (\ref{ur1}) by $u(x,t)$ and integrate
the resulting relation with respect to $x$ from  $0$ to $l$:

\begin{equation}
\int\limits_0^lu\partial_{0t}^\alpha
udx-\int\limits_0^lu(ku_x)_xdx+\int\limits_0^lqu^2dx=\int\limits_0^lufdx.
 \label{ur4}
\end{equation}

Let us transform the terms occurring in identity (\ref{ur4})

\begin{equation}
-\int\limits_0^lu(ku_x)_xdx=\int\limits_0^lku_x^2dx,\quad
\left|\int\limits_0^lufdx\right|\leq
\varepsilon\|u\|_0^2+\frac{1}{4\varepsilon}\|f\|_0^2, \quad
\varepsilon>0; \label{ur511}
\end{equation}

by virtue of inequality (\ref{ur5}) we obtain

\begin{equation}
\int\limits_0^lu(x,t)\partial_{0t}^\alpha
u(x,t)dx\geq\int\limits_0^l\frac{1}{2}\partial_{0t}^\alpha
u^2(x,t)dx=\frac{1}{2}\partial_{0t}^\alpha\|u\|_0^2. \label{ur512}
\end{equation}

Identity (\ref{ur4}), with regard of the above-performed
transformations, implies the inequality

\begin{equation}
\frac{1}{2}\partial_{0t}^\alpha\|u\|_0^2+c_1\|u_x\|_0^2\leq
\varepsilon\|u\|_0^2+\frac{1}{4\varepsilon}\|f\|_0^2.
 \label{ur7}
\end{equation}

By virtue of the inequality $\|u\|_0^2\leq (l^2/2)\|u_x\|_0^2$ for
$\varepsilon=c_1/l^2$ from (\ref{ur7}), we obtain the inequality

\begin{equation}
\partial_{0t}^\alpha\|u\|_0^2+c_1\|u_x\|_0^2\leq
\frac{l^2}{2c_1}\|f\|_0^2.
 \label{ur8}
\end{equation}

By applying the fractional differentiation operator
$D_{0t}^{-\alpha}$ to both sides of inequality (\ref{ur8}), we
obtain the estimate (\ref{ur9.1}) with constant $M=\max\{l^2/(2c_1
), 1\}/\min\{1,c_1\}$.

It follows from the a priori estimate (\ref{ur9.1}) that the
solution of problem (\ref{ur1})-(\ref{ur3}) is unique and
continuously depends on the input data.

%%%%%%%%%%%%%%%%%%%%%%%%%%%%%%%%%%%%%%%%%%%%%%%%%%%%%%%%%%%%%%%%%%%%%%%%
%%%%%%%%%%%%%%%%%%%%%%%%%%%%%%%%%%%%%%%%%%%%%%%%%%%%%%%%%%%%%%%%%%%%%%%%
%\section{ÒÐÅÒÜß ÊÐÀÅÂÀß ÇÀÄÀ×À ÄËß ÓÐÀÂÍÅÍÈß ÄÈÔÔÓÇÈÈ ÄÐÎÁÍÎÃÎ ÏÎÐßÄÊÀ}
%%%%%%%%%%%%%%%%%%%%%%%%%%%%%%%%%%%%%%%%%%%%%%%%%%%%%%%%%%%%%%%%%%%%%%%%
%%%%%%%%%%%%%%%%%%%%%%%%%%%%%%%%%%%%%%%%%%%%%%%%%%%%%%%%%%%%%%%%%%%%%%%%

{\bf 1.2. Third Boundary Value Problem.} In problem
(\ref{ur1})-(\ref{ur3}), we replace the boundary conditions
(\ref{ur2}) by the conditions

\begin{equation}
\begin{cases}
k(0,t)u_x(0,t)=\beta_1(t)u(0,t)-\mu_1(t),\\
-k(l,t)u_x(l,t)=\beta_2(t)u(l,t)-\mu_2(t),\quad 0\leq t\leq T.
\label{ur10}
\end{cases}
\end{equation}

In the rectangle $\bar Q_T$, consider the third boundary value
problem (\ref{ur1}), (\ref{ur3}), (\ref{ur10}).

To obtain a priori estimates for solutions of various nonstationary
problems, we use the well-known Gronwall-Bellman lemma \cite[p.
152]{ladij} , whose generalization is provided by the following
assertion.

{\bf Lemma 2.} Let a nonnegative absolutely continuous function
$y(t)$ satisfy the inequality
\begin{equation}
\partial_{0t}^\alpha y(t)\leq c_1y(t)+c_2(t), \quad 0<\alpha\leq1,
 \label{ur15}
\end{equation}
for almost all $t$ in $[0,T]$, where $c_1>0$ and $c_2(t)$ is an
integrable nonnegative function on $[0,T]$. Then
\begin{equation}
y(t)\leq
y(0)E_{\alpha}(c_1t^{\alpha})+\Gamma(\alpha)E_{\alpha,\alpha}(c_1t^{\alpha})D_{0t}^{-\alpha}c_2(t),
 \label{ur16}
\end{equation}
where $E_{\alpha}(z)=\sum_{n=0}^{\infty}z^n/\Gamma(\alpha n+1)$ and
 $E_{\alpha,\mu}(z)=\sum_{n=0}^{\infty}z^n/\Gamma(\alpha n+\mu)$ are the Mittag-Leffler
functions.

{\bf Proof.} Let $\partial_{0t}^\alpha y(t)- c_1y(t)=g(t)$, then
(e.g., see \cite[p. 17]{Pskhu})

\begin{equation}
y(t)=y(0)E_{\alpha}(c_1t^{\alpha})
+\int_0^t(t-\tau)^{\alpha-1}E_{\alpha,\alpha}(c_1(t-\tau)^{\alpha})g(\tau)d\tau\
 \label{ur17}
\end{equation}

By virtue of the inequality $g(t)\leq c_2(t)$, the positivity of the
Mittag-Leffler function $E_{\alpha,\alpha}(c_1(t-\tau)^{\alpha})$
for given parameters, and the growth of the function
$E_{\alpha,\alpha}(t)$, from (\ref{ur17}), we obtain the inequality

$$
y(t)\leq y(0)E_{\alpha}(c_1t^{\alpha}) +
\int\limits_0^t(t-\tau)^{\alpha-1}E_{\alpha,\alpha}(c_1(t-\tau)^{\alpha})c_2(\tau)d\tau\leq
$$

$$
\leq
y(0)E_{\alpha}(c_1t^{\alpha})+\Gamma(\alpha)E_{\alpha,\alpha}(c_1t^{\alpha})D_{0t}^{-\alpha}c_2(t).
$$

The proof of the lemma is complete.

{\bf Theorem 2.} If , in addition to the assumptions of Theorem 1,
$\beta_i(t), \mu_i(t)\in C[0,T]$, $|\beta_i(t)|\leq \beta$, for all
$t\in [0,T]$, $i=1,2$, then the solution $u(x,t)$ of problem
(\ref{ur1}),(\ref{ur3}),(\ref{ur10}) admits the a priori estimate
\begin{equation}
\|u\|_0^2+D_{0t}^{-\alpha}\|u_x\|_0^2\leq M
\left(D_{0t}^{-\alpha}\|f\|_0^2+D_{0t}^{-\alpha}\mu_1^2(t)+D_{0t}^{-\alpha}\mu_2^2(t)+\|u_0(x)\|_0^2\right).
 \label{ur9}
\end{equation}

{\bf Proof.} Just as in the proof of Theorem 1, we multiply Eq.
(\ref{ur1}) by $u(x,t)$ and integrate the resulting relation with
respect to $x$ from $0$ to $l$. By transforming the terms occurring
in identity (\ref{ur511}), we obtain relations (\ref{ur512}) and
(\ref{ur511}) with $\varepsilon=1/2$ and
$$
-\int\limits_0^lu(ku_x)_xdx=\beta_1(t)u^2(0,t)+\beta_2(t)u^2(l,t)-\mu_1(t)u(0,t)-\mu_2(t)u(l,t)+\int\limits_0^lku_x^2dx,
$$

Identity (\ref{ur511}), with regard of the above-performed
transformations, acquires the form
$$
\frac{1}{2}\partial_{0t}^\alpha\|u\|_0^2+c_1\|u_x\|_0^2\leq
$$

\begin{equation}
\leq
-\beta_1(t)u^2(0,t)-\beta_2(t)u^2(l,t)+\mu_1(t)u(0,t)+\mu_2(t)u(l,t)+
\frac{1}{2}\|u\|_0^2+\frac{1}{2}\|f\|_0^2.
 \label{ur12}
\end{equation}

By virtue of the inequalities
$$\mu_1(t)u(0,t)\leq\frac{1}{2}u^2(0,t)+\frac{1}{2}\mu_1^2(t),\quad
\mu_2(t)u(l,t)\leq\frac{1}{2}u^2(l,t)+\frac{1}{2}\mu_2^2(t),$$
$$u^2(0,t), u^2(l,t)\leq
\varepsilon\|u_x\|_0^2+(1/\varepsilon+1/l)\|u\|_0^2,\quad
\varepsilon>0,$$ from (\ref{ur12}) with
$\varepsilon=c_1/(4\beta+2)$, we obtain the inequality

\begin{equation}
\partial_{0t}^\alpha\|u\|_0^2+c_1\|u_x\|_0^2\leq
M_1\left(\|u\|_0^2+ \mu_1^2(t)+\mu_2^2(t)+\|f\|_0^2\right).
 \label{ur13}
\end{equation}

By applying the fractional differentiation operator
$D_{0t}^{-\alpha}$ to both sides of inequality (\ref{ur13}), we
obtain the inequality

$$ \|u\|_0^2+D_{0t}^{-\alpha}\|u_x\|_0^2\leq$$
\begin{equation}\leq
M_2\left(D_{0t}^{-\alpha}\|u\|_0^2+
D_{0t}^{-\alpha}\mu_1^2(t)+D_{0t}^{-\alpha}\mu_2^2(t)+D_{0t}^{-\alpha}\|f\|_0^2+\|u_0(x)\|_0^2\right).
 \label{ur14}
\end{equation}

By eliminating the second term from the left-hand side of inequality
(\ref{ur14}) and by using Lemma 2, where
$y(t)=D_{0t}^{-\alpha}\|u(x,t)\|_0^2$, $\partial_{0t}^\alpha
y(t)=\|u(x,t)\|_0^2$ and $y(0)=0$, we obtain the inequality

\begin{equation}
D_{0t}^{-\alpha}\|u\|_0^2\leq M_3\left(
D_{0t}^{-2\alpha}\mu_1^2(t)+D_{0t}^{-2\alpha}\mu_2^2(t)+D_{0t}^{-2\alpha}\|f\|_0^2
+\frac{t^\alpha}{\Gamma(\alpha+1)}\|u_0(x)\|_0^2\right),
 \label{ur18}
\end{equation}
where $M_3=\Gamma(\alpha)E_{\alpha,\alpha}(M_2T^{\alpha})$.

Since the inequality $D_{0t}^{-2\alpha}h(t)\leq
({t^\alpha\Gamma(\alpha)}/{\Gamma(2\alpha)})D_{0t}^{-\alpha}h(t)$
holds for every nonnegative integrable function $h(t)$ on $[0,T]$,
it follows from (\ref{ur14}) and (\ref{ur18}) that the a priori
estimate (\ref{ur9}) is true.

%%%%%%%%%%%%%%%%%%%%%%%%%%%%%%%%%%%%%%%%%%%%%%%%%%%%%%%%%%%%%%%%%%%%%%%%
%%%%%%%%%%%%%%%%%%%%%%%%%%%%%%%%%%%%%%%%%%%%%%%%%%%%%%%%%%%%%%%%%%%%%%%%
\section{BOUNDARY VALUE PROBLEMS FOR THE FRACTIONAL WAVE EQUATION}
%%%%%%%%%%%%%%%%%%%%%%%%%%%%%%%%%%%%%%%%%%%%%%%%%%%%%%%%%%%%%%%%%%%%%%%%
%%%%%%%%%%%%%%%%%%%%%%%%%%%%%%%%%%%%%%%%%%%%%%%%%%%%%%%%%%%%%%%%%%%%%%%%

{\bf 2.1. First Boundary Value Problem.} In the rectangle $\bar
Q_T=\{(x,t):0\leq x\leq l, 0\leq t\leq T\}$, consider the first
boundary value problem

\begin{equation}\label{ur19}
\partial_{0t}^{1+\alpha} u=\frac{\partial}{\partial
x}\left(k(x,t)\frac{\partial u}{\partial x}\right)-q(x,t)u+f(x,t)
,\quad 0<x<l,\quad 0<t\leq T,
\end{equation}

\begin{equation}
u(0,t)=0,\quad u(l,t)=0,\quad 0\leq t\leq T, \label{ur20}
\end{equation}

\begin{equation}
u(x,0)=u_0(x),\quad u_t(x,0)=u_1(x),\quad 0\leq x\leq l,
\label{ur21}
\end{equation}
where $\partial_{0t}^{1+\alpha}
u(x,t)=\int_0^tu_{\tau\tau}(x,\tau)(t-\tau)^{-\alpha}d\tau/\Gamma(1-\alpha)$
is the Caputo fractional derivative of order $1+\alpha$,
$0<\alpha<1$.

Throughout the following, we assume that there exists a solution
$u(x,t)\in C^{2,2}(\bar Q)$ of problem (\ref{ur19})-(\ref{ur21}).

{\bf Theorem 3.} If $k(x,t)\in C^{1,1}(\bar Q_T)$, $q(x,t)\in
C^{0,1}(\bar Q_T)$, $f(x,t)\in C(\bar Q_T)$, $0<c_1\leq k(x,t)\leq
c_2$, $0<m_1\leq q(x,t)\leq m_2$ and $|k_t(x,t)|, |q_t(x,t)|\leq
c_3$ everywhere on $\bar Q_T$, then the solution $u(x,t)$ of problem
(\ref{ur19})-(\ref{ur21}) admits the a priori estimate

\begin{equation}
D_{0t}^{\alpha-1}\|u\|_0^2+\|u\|_{W_2^1(0,l)}^2\leq M
\left(\int\limits_0^t\|f\|_0^2d\tau+\|u_1(x)\|_0^2+\|u_0(x)\|_{W_2^1(0,l)}^2\right),
 \label{ur25}
\end{equation}
where $\|u\|_{W_2^1(0,l)}^2=\|u\|_0^2+\|u_x\|_0^2$.

{\bf Proof.} Let us multiply Eq. (\ref{ur19})  by $u_t(x,t)$ and
integrate the resulting relation with respect to $x$ from $0$ to
$l$,

\begin{equation}
\int\limits_0^lu_t\partial_{0t}^{1+\alpha}
udx-\int\limits_0^lu_t(ku_x)_xdx+\int\limits_0^lquu_tdx=\int\limits_0^lu_tfdx.
 \label{ur22.1}
\end{equation}

Let us transform the terms occurring in identity (\ref{ur22.1})

\begin{equation}
\int\limits_0^lu_t\partial_{0t}^{1+\alpha}
udx=\int\limits_0^lu_t\partial_{0t}^{\alpha} u_tdx\geq
\frac{1}{2}\partial_{0t}^{\alpha}\|u_t\|_0^2,
 \label{ur22.2}
\end{equation}

$$
-\int\limits_0^lu_t(ku_x)_xdx=\frac{1}{2}\frac{\partial}{\partial
t}\int\limits_0^lku_x^2dx-\frac{1}{2}\int\limits_0^lk_tu_x^2dx,
$$

\begin{equation}
\int\limits_0^lquu_tdx=\frac{1}{2}\frac{\partial}{\partial
t}\int\limits_0^lqu^2dx-\frac{1}{2}\int\limits_0^lq_tu^2dx,\quad
\left|\int\limits_0^lu_tfdx\right|\leq
\frac{1}{2}\|u_t\|_0^2+\frac{1}{2}\|f\|_0^2.
 \label{ur22.3}
\end{equation}

By taking into account the performed transformations, from identity
(\ref{ur22.1}), we obtain the inequality

\begin{equation}
\partial_{0t}^\alpha\|u_t\|_0^2+\frac{\partial}{\partial
t}\int\limits_0^lku_x^2dx+\frac{\partial}{\partial
t}\int\limits_0^lqu^2dx\leq
M_4\left(\|u_t\|_0^2+\|u\|_{W_2^1(0,l)}^2+\|f\|_0^2\right);
 \label{ur23.1}
\end{equation}

by integrating this relation with respect to $\tau$ from $0$ to $t$,
we obtain the inequality

$$
D_{0t}^{\alpha-1}\|u_t\|_0^2+\|u\|_{W_2^1(0,l)}^2\leq
$$

\begin{equation}
\leq
M_5\left(\int\limits_0^t(\|u_{\tau}\|_0^2+\|u\|_{W_2^1(0,l)}^2)d\tau+
\int\limits_0^t\|f\|_0^2d\tau+\|u_1(x)\|_0^2+\|u_0(x)\|_{W_2^1(0,l)}^2\right).
 \label{ur24.0}
\end{equation}

By omitting the first term on the left-hand side in inequality
(\ref{ur24.0}) and by using the Gronwall-Bellman lemma \cite[p.
152]{ladij}, where $y(t)=\int_0^t\|u\|_{W_2^1(0,l)}^2d\tau$,
$y'(t)=\|u\|_{W_2^1(0,l)}^2$ and $y(0)=0$, we obtain
\begin{equation}
\int\limits_0^t\|u\|_{W_2^1(0,l)}^2d\tau\leq
M_6\left(\int\limits_0^t\left(\|u_{\tau}\|_0^2+
\|f\|_0^2\right)d\tau+\|u_1(x)\|_0^2+\|u_0(x)\|_{W_2^1(0,l)}^2\right).
 \label{ur24.1}
\end{equation}

Then, by omitting the second term on the left-hand side in
inequality (\ref{ur24.0}) and by using inequality (\ref{ur24.1}), we
obtain the inequality

\begin{equation}
D_{0t}^{\alpha-1}\|u_t\|_0^2\leq
M_7\left(\int\limits_0^t\|u_{\tau}\|_0^2d\tau+
\int\limits_0^t\|f\|_0^2d\tau+\|u_1(x)\|_0^2+\|u_0(x)\|_{W_2^1(0,l)}^2\right).
 \label{ur24.2}
\end{equation}

By Lemma 2, where $y(t)=\int_0^t\|u_{\tau}(x,\tau)\|_0^2d\tau$,
$\partial_{0t}^{\alpha}y(t)=D_{0t}^{\alpha-1}\|u_t(x,t)\|_0^2$ and
$y(0)=0$, from (\ref{ur24.2}) we obtain the inequality

\begin{equation}
\int\limits_0^t\|u_{\tau}\|_0^2d\tau\leq M_8\left(
D_{0t}^{-1-\alpha}\|f\|_0^2+\|u_1(x)\|_0^2+\|u_0(x)\|_{W_2^1(0,l)}^2\right).
 \label{ur24.3}
\end{equation}

By virtue of the inequality $D_{0t}^{-1-\alpha}\|f\|_0^2\leq
(t^{\alpha}/\Gamma(1+\alpha))\int_0^t\|f\|_0^2d\tau$, it follows
from inequalities (\ref{ur24.0}), (\ref{ur24.1}) and (\ref{ur24.3})
that the a priori estimate(\ref{ur25}) holds. The a priori estimate
(\ref{ur25}) implies that the solution of problem
(\ref{ur19})-(\ref{ur21}) exists and continuously depends on the
input data.

{\bf 2.2. Third Boundary Value Problem.} In the rectangle $\bar Q_T$
, consider the third boundary value problem (\ref{ur19}),
(\ref{ur21}), (\ref{ur10}).

{\bf Òåîðåìà 4.} If , in addition to the assumptions of Theorem 3,
$\beta_i(t), \mu_i(t)\in C^1[0,T]$, $\beta_i(t)\geq \beta>0$ and
$|\beta_{it}(t)|\leq c_4$ for all $t\in [0,T]$, $i=1,2$, then the
solution $u(x,t)$ of problem (\ref{ur19}), (\ref{ur21}),
(\ref{ur10}) admits the a priori estimate

$$
D_{0t}^{\alpha-1}\|u_t\|_0^2+\|u\|_{W_2^1(0,l)}^2\leq M\left(
\int\limits_0^t(\|f\|_0^2+\mu_{1\tau}^2(\tau)+\mu_{2\tau}^2(\tau))d\tau\right)+
$$

\begin{equation}
+M\left(\|\mu_1(t)\|_{C[0,T]}^2
+\|\mu_2(t)\|_{C[0,T]}^2+\|u_1(x)\|_0^2+\|u_0(x)\|_{W_2^1(0,l)}^2\right).
 \label{ur27}
\end{equation}

{\bf Proof.} Just as in Theorem 3, by multiplying Eq. (\ref{ur19})
by $u_t(x,t)$ and by integrating the resulting relation with respect
to $x$ from $0$ to $l$, we obtain identity (\ref{ur22.1}). By
transforming the terms occurring in identity (\ref{ur22.1}), we
obtain relations (\ref{ur22.2}) and (\ref{ur22.3}) and

$$
-\int\limits_0^lu_t(ku_x)_xdx=
$$

$$
=\frac{1}{2}\frac{\partial}{\partial
t}\left(\beta_1(t)u^2(0,t)+\beta_2(t)u(l,t)-2\mu_1(t)u(0,t)-2\mu_2(t)u(l,t)+\int\limits_0^lku_x^2dx\right)-
$$

$$
-\frac{1}{2}\beta_{1t}(t)u^2(0,t)-\frac{1}{2}\beta_{2t}(t)u^2(l,t)+\mu_{1t}(t)u(0,t)+\mu_{2t}(t)u(l,t)-\frac{1}{2}\int\limits_0^lk_tu_x^2dx.
$$

By taking into account the performed transformations and by using
the inequalities $$\mu_{1t}(t)u(0,t)\leq
\frac{1}{2}\mu_{1t}^2(t)+\frac{1}{2}u^2(0,t), \quad
\mu_{2t}(t)u(l,t)\leq
\frac{1}{2}\mu_{2t}^2(t)+\frac{1}{2}u^2(l,t),$$
$$u^2(0,t), u^2(l,t)\leq \|u_x\|_0^2+(1+1/l)\|u\|_0^2,$$ from
identity (\ref{ur22.1}), we obtain the inequality

$$
\partial_{0t}^\alpha\|u_t\|_0^2+\frac{\partial}{\partial
t}\left(\int\limits_0^lku_x^2dx+\int\limits_0^lqu^2dx\right)+
$$

$$
+\frac{\partial}{\partial
t}(\beta_1(t)u^2(0,t)+\beta_2(t)u(l,t)-2\mu_1(t)u(0,t)-2\mu_2(t)u(l,t))\leq
$$

\begin{equation}
\leq
M_9\left(\|u_t\|_0^2+\|u\|_{W_2^1(0,l)}^2+\|f\|_0^2+\mu_{1t}^2(t)+\mu_{2t}^2(t)\right).
 \label{ur23}
\end{equation}

By integrating inequality (\ref{ur23}) with respect to $\tau$ from
$0$ to $t$ and by taking into account the inequalities
$$2\mu_1(t)u(0,t)\leq \varepsilon
u^2(0,t)+(1/\varepsilon)\mu_1^2(t),\quad 2\mu_2(t)u(l,t)\leq
\varepsilon u^2(l,t)+(1/\varepsilon)\mu_2^2(t),$$ for
$\varepsilon=\beta$, $$u_0^2(0),
u_0^2(l)\leq(1+1/l)\|u_0(x)\|_{W_2^1(0,l)}^2,$$ we obtain

$$
D_{0t}^{\alpha-1}\|u_t\|_0^2+\|u\|_{W_2^1(0,l)}^2\leq
$$

$$
\leq
M_{10}\left(\int\limits_0^t(\|u_{\tau}\|_0^2+\|u\|_{W_2^1(0,l)}^2)d\tau+\int\limits_0^t(\|f\|_0^2+\mu_{1\tau}^2(\tau)+\mu_{2\tau}^2(\tau))d\tau\right)+
$$

\begin{equation}
+M_{10}\left(\|\mu_1\|_{C[0,T]}^2+\|\mu_2\|_{C[0,T]}^2+\|u_1(x)\|_0^2+\|u_0(x)\|_{W_2^1(0,l)}^2\right).
 \label{ur24}
\end{equation}

By analogy with the first boundary value problem, by using first the
Gronwall-Bellman lemma \cite[p. 152]{ladij} and then Lemma 2, from
inequality (\ref{ur24}), we obtain the a priori estimate
(\ref{ur27}).

%%%%%%%%%%%%%%%%%%%%%%%%%%%%%%%%%%%%%%%%%%%%%%%%%%%%%%%%%%%%%%%%%%%%%%%%%%%%%%%%%%%%%%%%%%%%%%%%
%%%%%%%%%%%%%%%%%%%%%%%%%%%%%%%%%%%%%%%%%%%%%%%%%%%%%%%%%%%%%%%%%%%%%%%%%%%%%%%%%%%%%%%%%%%%%%%%
%%%%%%%%%%%%%%%%%%%%%%%%%%%
%%%% Ñïèñîê ëèòåðàòóðû %%%%%
%%%%%%%%%%%%%%%%%%%%%%%%%%%%

\end{document}